\documentstyle[11pt]{article}
\def\arcsin{\mathop{\rm arcsin}\nolimits}
\def\arccos{\mathop{\rm arccos}\nolimits}
\def\arctan{\mathop{\rm arctan}\nolimits}
\def\arccot{\mathop{\rm arccot}\nolimits}
\def\arcsec{\mathop{\rm arcsec}\nolimits}
\def\arccsc{\mathop{\rm arccsc}\nolimits}

\def\erf{\mathop{\rm erf}\nolimits}
\newcommand{\one}[2]{\\[1mm] \noindent \begin{tabular}
{p{.5cm}p{12.5cm}}
\mbox{\rm (#1)}&
{$\displaystyle{#2}$}
\end{tabular}}

\newcommand{\C}{{\rm {\mbox{C{\llap{{\vrule
height1.52ex}\kern.4em}}}}}}
\newcommand{\N} {{\rm {\mbox{\protect\makebox[.15em][l]{I}N}}}}
\newcommand{\Q} {{\rm {\mbox{Q{\llap{{\vrule
height1.5ex}\kern.5em}}}}}}
\newcommand{\R} {{\rm {\mbox{\protect\makebox[.15em][l]{I}R}}}}

\newcommand{\Z} {{\rm {\mbox{\protect\makebox[.2em][l]{\sf Z}\sf
Z}}}}

\newcommand{\ed}[1]{\frac{1}{#1}}

\newcommand{\subs}[2]{\left. \makebox{\rule{0in}{2.5ex}} #1
\rb_{#2}}

\newcommand{\funkdeff}[4]{\left\{\begin{array}{ccc}
                                 #1 && \mbox{\rm{if} $#2$ } \\
                                 #3 && \mbox{\rm{if} $#4$ } 
                                 \end{array}
                          \right.}
\newcommand{\ueber}[2]{
                       \Big( \!
                       {{\small
                       \begin{array}{c}
                          #1\\
                          #2
                          \end{array}
                       }}
                       \! \Big) }

\newcommand{\pr}{\vspace{-2mm}\absatz{{\sl Proof:}}\hspace{5mm}}
\newcommand{\eop}{\hfill$\Box$\par\medskip\noindent}
\newcommand{\absatz}{\par\medskip\noindent}

\newcommand{\In}{\in}

\newcommand{\al}{\alpha}

\newcommand{\la}{\lambda}

\renewcommand{\phi}{\varphi}

\newcommand{\menge}[3]{\left\{#1 \In #2 \; \lb \; #3 \right.
\right\} }

\newcommand{\1}{{\bf{1}}}
\newcommand{\2}{{\bf{2}}}
\newcommand{\3}{{\bf{3}}}

\newcommand{\6}{{\bf{6}}}
\newcommand{\7}{{\bf{7}}}

\newcommand{\pf}{\rightarrow}

\newcommand{\be}{\begin{equation}}
\newcommand{\ee}{\end{equation}}
\newcommand{\bea}{\begin{eqnarray}}
\newcommand{\eea}{\end{eqnarray}}
\newcommand{\beao}{\begin{eqnarray*}}
\newcommand{\eeao}{\end{eqnarray*}}

\newcommand{\lk}{\left(}
\newcommand{\rk}{\right)}
\newcommand{\lb}{\left|}
\newcommand{\rb}{\right|}
\newcommand{\bi}{\bibitem}

\newcommand{\bT}{\begin{theorem}}
\newcommand{\eT}{\end{theorem}}
\newcommand{\bL}{\begin{lemma}}
\newcommand{\eL}{\end{lemma}}
\newcommand{\bC}{\begin{corollary}}  
\newcommand{\eC}{\end{corollary}}
\newcommand{\bt}{\begin{tabbing} 12345 \= \kill}
\newcommand{\et}{\end{tabbing}}

\newcommand{\abs}{\\[3mm]}

\newtheorem{theorem}{Theorem}
\newtheorem{algorithm}{Algorithm}
\newtheorem{lemma}{Lemma}
\newtheorem{corollary}{Corollary}
\newtheorem{remark}{Remark}

\begin{document}
\begin{center}{
{\LARGE {\bf {Spaces of Functions Satisfying}}}
\vspace{3mm}\\
{\LARGE {\bf {Simple Differential Equations}}}
\vspace{5mm}\\
{\large {\sc {Wolfram Koepf}}}
\vspace{3mm}\\
{\sl Konrad-Zuse-Zentrum f\"ur Informationstechnik Berlin,
Heilbronner Str. 10,
D-10711 Berlin,}
\\
{\sl Federal Republic of Germany}
\vspace{5mm}\\
{\large {\sc {Dieter Schmersau}}}
\vspace{3mm}\\
{\sl Freie Universit\"at Berlin, Fachbereich Mathematik, Arnimallee
3,
D-14195 Berlin,}
\\
{\sl Federal Republic of Germany}
}
\\[3mm]
Konrad-Zuse-Zentrum Berlin (ZIB), Technical Report TR 94-2, 1994
\end{center}
\vspace{0.5cm}
\begin{center}
{\bf {Abstract:}}
\end{center}
\begin{center}
\begin{tabular}{p{13cm}}
{{\small
In \cite{Koe92}--\cite{Koe93c} the first author published an
algorithm for the
conversion of analytic functions for which derivative rules are
given
into their representing 
power series $\sum\limits_{k=0}^{\infty}a_{k}z^{k}$ at the origin
and
vice versa, implementations of which exist in {\sc Mathematica}
\cite{Wol}, (s.\ \cite{Koe93c}), {\sc Maple} \cite{Map} (s.\
\cite{GK})
and {\sc Reduce} \cite{Red} (s.\ \cite{Neun}).

One main part of this procedure is an algorithm to derive a
homogeneous
linear differential equation with polynomial coefficients
for the given function. We call this type of ordinary differential
equations {\sl simple}.

Whereas the opposite question to find functions satisfying given
differential equations is studied in great detail, our question
to find differential equations that are satisfied by given
functions
seems to be rarely posed.

In this paper we consider the family $F$ of functions satisfying
a simple differential equation generated by the rational, the
algebraic,
and certain transcendental functions. It turns out that $F$ forms
a 
linear space of transcendental functions. 
Further $F$ is closed under multiplication and under the
composition with rational functions and rational powers.
These results had been published by
Stanley who had proved them by theoretical algebraic
considerations.

In contrast our treatment is purely algorithmically oriented.
We present algorithms that generate simple differential equation
for $f+g$, $f\cdot g$, $f\circ r$ ($r$ rational), and $f\circ
x^{p/q}$ 
($p,q\in\N_0$), given simple differential equations for $f$, and
$g$, and
give a priori estimates for the order of the resulting differential
equations.
We show that all order estimates are sharp.

After finishing this article we realized
that in independent work Salvy and Zimmermann published similar
algorithms.
Our treatment gives a detailed description of those algorithms and
their
validity.
}}
\end{tabular}
\\[1cm]
\end{center}

\section{Simple functions}

Many mathematical functions satisfy a homogeneous
linear differential equation with polynomial coefficients.
We call such an ordinary differential equation {\sl simple}. Also
a function $f$ that satisfies a simple differential equation 
is called simple.
The least order of such a simple differential equation fulfilled by
$f$
is called the order of $f$. 

Technically there is no difference between a differential equation
with
polynomial coefficients, and one with rational coefficients, as we
can
multiply such a differential equation by its common denominator, so
we
may consider the case that the coefficients are members of the
field
$K[x]$ where $K[x]$ is one of $\Q[x]$, $\R[x]$, or $\C[x]$.

Examples of simple functions are
\begin{enumerate}
\item
all rational functions $p/q$ ($p,q$ polynomials) are simple of
order $1$:
$p\,q\,f'+(p\,q'-q\,p')\,f=0$,
\item
all algebraic functions $f$ are simple
(see e.\ g.\ \cite{CC1}--\cite{CC2},\cite{Wal1}, and \cite{Koe922})
of the order of $f$,
\item
the power function $x^\alpha$ satisfies the simple differential
equation
$x\,f'-\alpha\,f=0$,
\item
the exponential function $e^x$ satisfies the simple differential
equation $f'-f=0$,
\item
the logarithm function $\ln x$ satisfies the simple differential
equation $x\,f''+f'=0$,
\item
the sine function $\sin x$ and the cosine function $\cos x$
satisfy the simple differential equation $f''+f=0$,
\item
the inverse sine function $\arcsin x$ and the inverse cosine
function
$\arccos x$ satisfy the simple differential equation
$(x^2-1)\,f''+x\,f'=0$,
\item
the inverse tangent function $\arctan x$ and the inverse cotangent
functions $\arccot x$ satisfy the simple differential equation
$(1+x^2)\,f''+2\,x\,f'=0$,
\item
the inverse secant function $\arcsec x$ and the inverse cosecant
function $\arccsc x$ satisfy the simple differential equation
$(x^3-x)\,f''+(2x^2-1)\,f'=0$,
\item
the error function $\erf x$ satisfies the simple differential
equation
$f''+2\,x\,f'=0$.
\end{enumerate}
Note that for rational functions with rational coefficients, 
algebraic functions given by a rational coefficient equation, 
the power function for $\al \in\Q$, and in all other examples
we have $K[x]=\Q[x]$. 

Also many special functions are simple with $K[x]=\Q[x]$:
Airy functions (see e.\ g.\ \cite{AS}, \S~10.4), Bessel functions
(see e.\
g.\ \cite{AS}, Ch.~9--11), all kinds of orthogonal polynomials
(see e.\ g.\ \cite{Sze}, \cite{Tri}),
and functions of hypergeometric type (see \cite{Koe92}).

Not all elementary transcendental functions, however, are simple,
the simplest example of which probably is the tangent function
$\tan x$,
compare \cite{Sta}, Example 2.5.
\bT
\label{th:tanx}
{\rm
The functions $\tan x$ and $\sec x$ do not satisfy simple
differential 
equations.
}
\eT
\pr
It is easily seen that the tangent function $f(x)=\tan x$
satisfies the nonlinear differential equation
\be
f'=1+f^2
\;.
\label{eq:nonlineartangent}
\ee
Differentiation of (\ref{eq:nonlineartangent}) yields after further
substitution of (\ref{eq:nonlineartangent})
\[
f''=(1+f^2)'=2\,f\,f'=2\,f\,(1+f^2)
\;,
\]
and inductively we get representations of
$f^{(k)}=P_k(f)\;(k\in\N)$ by
polynomial expressions $P_k$ in $f$.

Assume now, the tangent function would satisfy the simple
differential
equation
\be
\sum\limits_{k=0}^n
p_k\,f^{(k)}=0\quad\quad(p_k\;\;\mbox{polynomials})
\;,
\label{eq:simple DE}
\ee
then we can substitute $f^{(k)}$ by $P_k$ arriving at an algebraic
identity for the tangent function. This obviously is a
contradiction to
the fact that the tangent function is transcendental, and therefore
a
differential equation of type (\ref{eq:simple DE}) cannot be
satisfied.

Let us now consider $f(x)=\sec x$. In this case we may establish
the nonlinear
differential equation
\[
f'^2=f^4-f^2
\]
from which, by differentiation, the second order differential
equation
\be
f''=2f^3-f
\label{eq:secsecondorder}
\ee
follows. If now a simple differential equation (\ref{eq:simple DE})
would be valid for $f$, then also
\[
p_1\,f'=-\sum_{\begin{array}{c}
\\[-6.5mm]
{\scriptstyle{k=0}}\\[-2.5mm]
{\scriptstyle{k\neq 1}}
\end{array}}^n
p_k\,f^{(k)}
\;,
\]
and taking the square we get
\[
p_1^2\,f'^2=p_1^2\,(f^4-f^2)=\lk
\sum_{\begin{array}{c}
\\[-6.5mm]
{\scriptstyle{k=0}}\\[-2.5mm]
{\scriptstyle{k\neq 1}}
\end{array}}^n
p_k\,f^{(k)}
\rk^2
\;.
\]
If we substitute (\ref{eq:secsecondorder}), and the corresponding 
representations for $f^{(k)}\;(k=3,\ldots,n)$ which we get by
differentiating
(\ref{eq:secsecondorder}), we arrive at
an algebraic identity, and the conclusion follows by the
transcendency of $f$, again.
\eop
By similar means we see that the elementary functions $\cot x$, 
and $\csc x$ as well as the corresponding hyperbolic functions are
not
simple. Moreover, the generating functions
\[
\frac{x}{e^x-1}=\sum_{k=0}^\infty B_k\,\frac{x^k}{k!}\;,
\quad\quad\quad\mbox{and}\quad\quad\quad
\frac{e^{x/2}}{e^x+1}=\sum_{k=0}^\infty
\frac{E_k}{2^{k+1}}\,\frac{x^k}{k!}
\]
of the Bernoulli and Euler numbers $B_k$, and $E_k\;(k\in\N_0)$
(see e.\ g.\
\cite{AS}, (23.1)) satisfy the nonlinear differential equations
\[
f'=\ed x\lk (1-x)\,f-f^2\rk,
\quad\quad\quad\mbox{and}\quad\quad\quad
f'^2=\frac{f^2}{4}-f^4
\;,
\]
respectively, and so by a similar procedure are realized not to be
simple. 
As exactly the simple differential equations correspond to
simple recurrence equations for the Taylor series coefficients (see
e.\ g.\
\cite{Koe92}), this implies the following
\bC
{\rm
Neither the Bernoulli nor the Euler numbers satisfy a finite
homogeneous linear
recurrence equation with polynomial coefficients.
}
\eC

\section{Construction of simple functions}

If simple functions are given, then we may construct further
examples of 
simple functions by the following procedures: integration,
differentiation, addition, multiplication, and the composition with
certain functions, namely rational functions and rational powers.
Before we prove this main result, we consider an important
particular case. We will prove that if $f$ is simple, and $q\in\N$,
then 
$f\circ x^{1/q}$ is simple, too. To get this result, we need the
following
\bL
\label{l:1/q}
{\rm
Let $f:D\pf\R$ be infinitely often differentiable, let be
$D\subset\R^+\neq
\emptyset$, and let $q\in\N$ with $q>2$ be given. Consider the
function
$h:D\pf\R$ with $h(x):=f(x^{1/q})$. 
\begin{enumerate}
\item[(a)]
Then we have for $r\in\N$ with $1\leq r \leq q-1$
\be
x^{r/q}\,f^{(r)}(x^{1/q})=\sum_{k=1}^r C_{qrk}^0\,x^k\,h^{(k)}(x)
\label{eq:thetan}
\ee
with constants $C_{qrk}^0$, in particular
\[
C_{q11}^0=q\;,
\quad\quad
C_{qr1}^0=q\,(q-1)\cdots (q-r+1)\;,
\quad\quad\mbox{and}\quad\quad
C_{qrr}^0=q^r
\;,
\]
and
\be
f^{(q)}(x^{1/q})=\sum_{k=1}^q C_{q0k}^1\,x^{k-1}\,h^{(k)}(x)
\label{eq:scharf1/q}
\ee
with constants $C_{q0k}^1$, in particular
\[
C_{q01}^1=q!
\;,
\quad\quad\mbox{and}\quad\quad
C_{q0q}^1=q^q
\;.
\]
\item[(b)]
Further for $s\in\N$ we have
\[
f^{(sq)}(x^{1/q})=\sum_{k=s}^{sq} C_{q0k}^s\,x^{k-s}\,h^{(k)}(x)
\]
with constants $C_{q0k}^s$, in particular
\be
C_{q,0,sq}^s=q^{sq}
\;,
\label{eq:sq}
\ee
and for $r\in\N$ with $1\leq r\leq q-1$
\[
x^{r/q}\,f^{(sq+r)}(x^{1/q})=\sum_{k=s+1}^{sq+r} 
C_{qrk}^s\,x^{k-s}\,h^{(k)}(x)
\]
with constants $C_{qrk}^s$, in particular
\be
C_{q,r,sq+r}^s=q^{sq+r}
\;.
\label{eq:sq+r}
\ee
\end{enumerate}
}
\eL
\pr
Differentiating the identity $f(x^{1/q})=h(x)$, we get
\[
x^{1/q}\,f'(x^{1/q})=q\,x\,h'(x)
\;,
\]
a further differentiation yields
\[
x^{2/q}\,f''(x^{1/q})+x^{1/q}\,f'(x^{1/q})=q^2\,x\,h'(x)+q^2\,x^2
\,h''(x)
\;,
\]
and therefore
\[
x^{2/q}\,f''(x^{1/q})=q\,(q-1)\,x\,h'(x)+q^2\,x^2\,h''(x)
\;.
\]
Differentiating once more and rearranging similarly, gives
\[
x^{3/q}\,f'''(x^{1/q})=q\,(q-1)\,(q-2)\,x\,h'(x)+3\,q^2\,(q-1)\,x
^2\,h''(x)+
q^3\,x^3\,h'''(x)
\;.
\]
Assume now, for $r\in\N$ with $1\leq r\leq q-1$ the relation
\[
x^{r/q}\,f^{(r)}(x^{1/q})=\sum_{k=1}^r C_{qrk}^0\,x^k\,h^{(k)}(x)
\]
with
\[
C_{qr1}^0=q\,(q-1)\cdots (q-r+1)\;,
\quad\quad\mbox{and}\quad\quad
C_{qrr}^0=q^r
\]
is valid, then by another differentiation we get
\beao
x^{(r+1)/q}\,f^{(r+1)}(x^{1/q})
&=&
(q-r)\,C_{qr1}^0\,x\,h'(x)
\\&=&
+\sum_{k=2}^r
\lk (k\,q-r)\,C_{qrk}^0+q\,C_{q,r,k-1}^0\rk x^k\,h^{(k)}(x)
+q\,C_{qrr}^0\,x^{r+1}\,h^{(r+1)}(x)
\;.
\eeao
For $r\leq q-1$, we set
\[
C_{q,r+1,1}^0:=(q-r)\,C_{qr1}^0=q\,(q-1)\cdots(q-r)
\;,
\]
\[
C_{q,r+1,r+1}^0:=q\,C_{qrr}^0=q\,q^r=q^{r+1}
\;,
\]
and
\[
C_{q,r+1,k}^0:=(k\,q-r)\,C_{qrk}^0+q\,C_{q,r,k-1}^0
\quad\quad(2\leq k\leq r)
\;.
\]
For $r=q$ the same method applies, if we set
\[
C_{q01}^1:=C_{q,q-1,1}^0=q!
\;,
\]
\[
C_{q0k}^1:=q\,C_{q,q-1,q-1}^0=q\,q^{q-1}=q^{q}
\;,
\]
and
\[
C_{q0k}^1:=((k-1)\,q+1)\,C_{q,q-1,k}^0+q\,C_{q,q-1,k-1}^0
\quad\quad(2\leq k\leq q-1)
\;,
\]
which proves (a). 

To prove (b) we assume that for some $s\in\N$
\[
f^{(sq)}(x^{1/q})=\sum_{k=s}^{sq} C_{q0k}^s\,x^{k-s}\,h^{(k)}(x)
\]
with
\[
C_{q,0,sq}^s=q^{sq}
\;.
\]
If we define the functions $f_s$ and $h_s$ by
\[
f_s(x):=f^{(sq)}(x)
\;,
\quad\quad\mbox{and}\quad\quad
h_s(x):=\sum_{k=s}^{sq} C_{q0k}^s\,x^{k-s}\,h^{(k)}(x)
\;,
\]
then we may apply (a), and get for $1\leq r\leq q-1$ that
\[
x^{r/q}\,f_s^{(r)}(x^{1/q})=\sum_{k=1}^r
C_{qrk}^0\,x^k\,h_s^{(k)}(x)
\;,
\]
and further
\[
f_s^{(q)}(x^{1/q})=\sum_{k=1}^q C_{q0k}^1\,x^{k-1}\,h_s^{(k)}(x)
\;,
\]
hence
\be
x^{r/q}\,f^{(sq+r)}(x^{1/q})=\sum_{k=1}^r
C_{qrk}^0\,x^k\,h_s^{(k)}(x)
\,,
\label{eq:stern1}
\ee
and
\be
f^{((s+1)q)}(x^{1/q})=\sum_{k=1}^q C_{q0k}^1\,x^{k-1}\,h_s^{(k)}(x)
\;.
\label{eq:stern2}
\ee
To finish our proof, we will calculate the derivatives 
$h_s^{(k)}\;(k=1,\ldots,q)$ algorithmically. We have
\[
h_s(x)=\sum_{j=s}^{sq} C_{q0j}^s\,x^{j-s}\,h^{(j)}(x)
\;,
\]
and therefore
\beao
h_s'(x)
&=&
\sum_{j=s}^{sq} C_{q0j}^s\,x^{j-s}\,h^{(j+1)}(x)
+
\sum_{j=s+1}^{sq} (j-s)\,C_{q0j}^s\,x^{j-s-1}\,h^{(j)}(x)
\\
&=&
\sum_{j=s+1}^{sq+1} C_{q,0,j-1}^s\,x^{j-s-1}\,h^{(j)}(x)
+
\sum_{j=s+1}^{sq} (j-s)\,C_{q0j}^s\,x^{j-s-1}\,h^{(j)}(x)
\\
&=&
\sum_{j=s+1}^{sq} \Big( (j-s)\,C_{q0j}^s+C_{q,0,j-1}^s \Big)
\,x^{j-s-1}\,h^{(j)}(x)
+
C_{q,0,sq}^s\,x^{sq-q}\,h^{(sq+1)}(x)
\;.
\eeao
We may write this as
\[
h_s'(x)
=
\sum_{j=s+1}^{sq+1} D_j^1\,x^{j-s-1}\,h^{(j)}(x)
\]
with
\[
D_j^1:=\funkdeff{(j-s)\,C_{q0j}^s+C_{q,0,j-1}^s}
{s+1\leq j\leq s\,q}
{C_{q,0,sq}^s=q^{sq}}{j=s\,q+1}
\;.
\]
Continuing iteratively in the same fashion, we get representations
\[
h_s^{(k)}(x)=\sum_{j=s+k}^{sq+k} D_j^k\,x^{j-s-k}\,h^{(j)}(x)
\]
with
\be
D_j^k=\funkdeff{D_{j-1}^{k-1} +(j-s-k+1)\,D_{j}^{k-1}}
{s+k\leq j\leq s\,q+k-1}
{q^{sq}}{j=s\,q+k}
\;.
\label{eq:plus}
\ee
We substitute these equations in (\ref{eq:stern1}) to get
\[
x^{r/q}\,f^{(sq+r)}(x^{1/q})=
\sum_{k=1}^r \sum_{j=s+k}^{sq+k}
C_{qrk}^0\,D_j^k\,x^{j-s}\,h^{(j)}(x)
\;.
\]
For $j$ with $s+1\leq j\leq s\,q+r$ we use the notation
\[
M_j:=\menge{k}{\N}{1\leq k\leq r\quad\mbox{and}\quad s+k\leq j\leq
s\,q+k}
\;.
\]
Then $M_{sq+r}=\{r\}$, and
\[
x^{r/q}\,f^{(sq+r)}(x^{1/q})=\sum_{j=s+1}^{sq+r} 
C_{qrj}^s\,x^{j-s}\,h^{(j)}(x)
\]
with
\[
C_{qrj}^s:=\sum_{k\in M_j} C_{qrk}^0\,D_j^k
\]
for $s+1\leq j\leq s\,q+r$, in particular
\[
C_{q,r,sq+r}^s=C_{qrr}^0\,D_{sq+r}^r=q^r\,q^{sq}=q^{sq+r}
\]
proving (\ref{eq:sq+r}).

To prove (\ref{eq:sq}), we substitute Equations (\ref{eq:plus})
in (\ref{eq:stern2}) to get
\[
f^{((s+1)q)}(x^{1/q})=
\sum_{j=s+1}^{(s+1)q} C_{q0j}^{s+1}\,x^{j-s-1}\,h^{(j)}(x)
\]
with
\[
C_{q0j}^{s+1}:=\sum_{k\in M_j} C_{q0k}^1\,D_j^k
\]
for $s+1\leq j\leq (s+1)\,q$, in particular
\[
C_{q,0,(s+1)q}^{s+1}=C_{q0q}^1\,D_{sq+q}^q=q\,q^{sq}=q^{(s+1)q}
\]
finishing our proof.
\eop
\begin{remark}
{\rm
We remark that one easily can establish that all
constants $C_{qrk}^s$ are integers. Further, we note that 
using the differential operator $\theta_q:=q\,x\,\frac{d}{dx}$
(compare with the proof of \cite{Koe92}, Theorem 8.1) Equation
(\ref{eq:thetan})
gets the simple form
\[
x^{r/q}\,f^{(r)}(x^{1/q})=\theta_q\,(\theta_q-1)\,(\theta_q-2)\cdots
\,(\theta_q-(r-1))\,h(x)
\quad\quad
(1\leq r\leq q)
\;.
\]
}
\end{remark}

Next we give some applications of the lemma.
First assume that $f'=f$, a solution of which is given by the
exponential
function $f(x)=e^x$. Then $f^{(q)}=f$, and so
$f^{(q)}(x^{1/q})=f(x^{1/q})$,
and therefore we get for $h(x):=e^{x^{1/q}}$ the differential
equation
\[
\sum_{k=1}^q C_{q0k}^1\,x^{k-1}\,h^{(k)}=h
\;,
\]
e.\ g.\ for $q=3$ we have the differential equation for
$e^{x^{1/3}}$
\[
27\,x^2\,h'''+54\,x\,h''+6\,h'=h
\;.
\]
Further we consider the homogeneous Euler differential equation
\be
\sum_{k=0}^n a_k\,x^k\,f^{(k)}=0
\label{eq:Euler differential equation}
\ee
with constants $a_k\;(k=0,\ldots,n)$, $a_n\neq 0$. As by the Lemma
we have
\beao
\subs{t^r\,f^{(r)}(t)}{t=x^{1/q}}
&=&
\sum_{k=1}^r C_{qrk}^0\,x^k\,h^{(k)}(x)
\quad(1\leq r \leq q-1)
\;,
\\[2mm]
\subs{t^q\,f^{(q)}(t)}{t=x^{1/q}}
&=&
\sum_{k=1}^q C_{q0k}^1\,x^k\,h^{(k)}(x)
\;,
\\[2mm]
\subs{t^{sq}\,f^{(sq)}(t)}{t=x^{1/q}}
&=&
\sum_{k=s}^{sq} C_{q0k}^s\,x^k\,h^{(k)}(x)
\quad(s\in\N)
\;,
\eeao
and
\[
\subs{t^{sq+r}\,f^{(sq+r)}(t)}{t=x^{1/q}}=
\sum_{k=s+1}^{sq+r} C_{qrk}^s\,x^k\,h^{(k)}(x)
\quad(1\leq r \leq q-1)
\;,
\quad
\]
we immediately conclude the following
\bC
\label{cor:euler}
{\rm
Assume, the function $f$ satisfies a homogeneous Euler differential
equation
(\ref{eq:Euler differential equation}) of order $n$. Then the
function
$f\circ x^{1/q}$ satisfies a homogeneous Euler differential
equation of
the same order $n$.
\hfill$\Box$
}
\eC
As an application of Corollary~\ref{cor:euler} we get
\bC
{\rm
Assume $P$ is a polynomial of degree $n$, then $P\circ x^{1/q}$ 
satisfies a homogeneous Euler differential equation of order $\leq
n+1$.
}
\eC
\pr
$P$ satisfies the differential equation $x^{n+1}\,f^{(n+1)}=0$.
\eop
Our next step is to prove another lemma which we will need to prove
that if 
$f$ is simple, and $q\in\N$, then so is $f\circ x^{1/q}$.
\bL
\label{l:lin_nq}
{\rm
Let $f:D\pf\R$ be an infinitely often differentiable function that
is simple
of order $n$, let be $D\subset\R^+\neq
\emptyset$, and let $q\in\N$ with $q>2$ be given. Consider the
function
$h:D\pf\R$ with $h(x):=f(x^{1/q})$. Then the linear space
$L_f$ over $K[x]$ (i.e.\ with rational function coefficients)
generated by the
functions $x^{r/q}\,f^{(k)}(x^{1/q})\;(r,k\in\N_0)$ has dimension 
$\leq n\,q$.
}
\eL
\pr
Given simple $f$ of order $n$, the linear space $L_f$ is generated
by the functions of the $\N_0\times\N_0$ matrix
\[
\begin{array}{cccccc}
f(y),&y\,f(y),&y^2\,f(y),&\ldots,&y^r\,f(y),&\ldots\\[1mm]
f'(y),&y\,f'(y),&y^2\,f'(y),&\ldots,&y^r\,f'(y),&\ldots\\[1mm]
f''(y),&y\,f''(y),&y^2\,f''(y),&\ldots,&y^r\,f''(y),&\ldots\\[1mm]
\vdots&\vdots&\vdots&&\vdots&\\[1mm]
f^{(k)}(y),&y\,f^{(k)}(y),&y^2\,f^{(k)}(y),&\ldots,&y^r\,f^{(k)}(
y),&\ldots
\\[1mm]
\vdots&\vdots&\vdots&&\vdots
\end{array}
\]
using the abbreviation $y:=x^{1/q}$.

At most the first $q$ expressions 
$f(x^{1/q}),x^{1/q}\,f(x^{1/q}),\ldots,x^{(q-1)/q}\,f(x^{1/q})$
of the first row
are linearly independent (over $K[x]$) as the expressions for
$r=sq\;(s\in\N)$
depend linearly on the first expression $f(y)$ ($y^{sq}=x^s$ is a
polynomial),
further the expressions for $r=sq+1\;(s\in\N)$
depend linearly  on the second expression $y\,f(y)$, and so on.
The same conclusion follows for each other row so that we obtain
that $L_f$
is generated by the $q\times \N_0$ matrix $M$
\[
\begin{array}{ccccc}
f(y),&y\,f(y),&y^2\,f(y),&\ldots,&y^{q-1}\,f(y),\\[1mm]
f'(y),&y\,f'(y),&y^2\,f'(y),&\ldots,&y^{q-1}\,f'(y),\\[1mm]
f''(y),&y\,f''(y),&y^2\,f''(y),&\ldots,&y^{q-1}\,f''(y),\\[1mm]
\vdots&\vdots&\vdots&&\vdots\\[1mm]
f^{(k)}(y),&y\,f^{(k)}(y),&y^2\,f^{(k)}(y),&\ldots,&y^{q-1}\,f^{(
k)}(y),
\\[1mm]
\vdots&\vdots&\vdots&&\vdots
\end{array}
\;.
\]
Now, we consider the first column. We will show that at most the
first
$n$ expressions of this column are linearly independent. 
By the given simple differential equation for $f$ there are
polynomials $p$, 
and $p_k\;(k=0,\ldots,n-1)$ such that
\[
p(y)\,f^{(n)}(y)=\sum_{k=0}^{n-1} p_k(y)\,f^{(k)}(y)
\;,
\]
or
\be
p(x^{1/q})\,f^{(n)}(x^{1/q})=\sum_{k=0}^{n-1}
p_k(x^{1/q})\,f^{(k)}(x^{1/q})
\;.
\label{eq:vorerst}
\ee
If $p(x^{1/q})$ turns out to be a polynomial in $x$, then this
equation
tells that $f^{(n)}(x^{1/q})$ depends linearly (over $K[x]$) on the
elements $x^{r/q}\,f^{(k)}(x^{1/q})\;(r\in\N_0, \;0\leq k\leq n-1)$
of the rows above in our scheme. As this, however, is not
necessarily the case,
we construct a polynomial $r$ such that $p(x^{1/q})\,r(x^{1/q})$ is
a 
polynomial in $x$.

Therefore we use the complex factorization 
$p(y)=c\,(y-y_1)\,(y-y_2)\cdots (y-y_m)$ of $p$, and by the
identity
\[
y^q-y_l^q=(y-y_l)\,\sum_{j=0}^{q-1} y^j\,y_l^{q-1-j}
\]
we see that multiplication of $p$ by the polynomial
\[
r(y):=\prod_{l=1}^m\lk\sum_{j=0}^{q-1} y^j\,y_l^{q-1-j}\rk
\]
yields a polynomial $s$ in the variable $x^q$. Even though, for
technical
reasons, we use a complex factorization, it turns out that 
$r\in \R[x] \;(\Q[x])$ if $p\in \R[x] \;(\Q[x])$. Hence
$P(x):=s(x^{1/q})=p(x^{1/q})\,r(x^{1/q})$ is a polynomial in $x$.
Thus
we see that multiplying (\ref{eq:vorerst}) by $r(y)$ generates the 
representation
\[
P(x)\,f^{(n)}(x^{1/q})=\sum_{k=0}^{n-1}
P_k(x^{1/q})\,f^{(k)}(x^{1/q})
\]
with polynomials $P$, and $P_k\;(k=0,\ldots,n-1)$.

This equation tells that $f^{(n)}(x^{1/q})$ depends linearly on the
functions of the rows above. Similarly, by an induction argument,
$f^{(j)}(x^{1/q})\;(j>n)$ depend linearly on those functions
$x^{r/q}\,f^{(k)}(x^{1/q})$ $(0\leq r\leq q-1,0\leq k\leq n-1)$
of the first $n$ rows of $M$.

Considering the other columns, the same argument can be applied to
show
that $y^r\,f^{(j)}(y)$ depend linearly on the functions of the
first $n$ 
rows of $M$. This shows finally that $L_f$ is generated by the
$n\,q$ functions
$x^{r/q}\,f^{(k)}(x^{1/q})\;(0\leq r\leq q-1,0\leq k\leq n-1)$, and
therefore has dimension $\leq n\,q$.
\eop
Lemma \ref{l:1/q} shows immediately
that for $h:=f\circ x^{1/q}$ the linear space $L_h$ over $K[x]$
generated
by $h, h', h'',\ldots$ is a subset of $L_f$ (declared in Lemma
\ref{l:lin_nq}),
and so by Lemma \ref{l:lin_nq} is of dimension $\leq n\,q$. Thus we
have proved
the following Theorem (compare \cite{Sta}, Theorem 2.7, \cite{SZ},
{\sc Maple}
function {\tt algebraicsubs}).
\bT
\label{th:1/q}
{\rm
Let $f$ be simple of order $n$, and let $q\in\N$. Then 
$f\circ x^{1/q}$ is simple of order $\leq n\,q$.
\hfill$\Box$
}
\eT
Whereas Theorem~\ref{th:1/q} gives a complete answer with regard to
the
existence of a simple differential equation, in \cite{Sta}, Theorem
2.7, it
is shown that this is true for the composition with each algebraic
function
$\phi(x)$ with $\phi(0)=)$), it does not tell anything
how such a differential equation may be obtained. Therefore we
continue to
describe an algorithm that generates the differential equation of
lowest 
order for $h$, which gives a second proof for Theorem~\ref{th:1/q}.
\abs
{\bf Algorithm to Theorem \ref{th:1/q}}:
Given the differential equation
\be
\sum_{k=0}^n p_k(x)\,f^{(k)}(x)=0
\quad\quad\quad(p_k\;\mbox{polynomials}\;(k=0,\ldots,n),\quad
p_n\not\equiv 0)
\label{eq:givenDE}
\ee
for $f$ we have to construct a simple differential equation for 
$h(x):=f(x^{1/q})\;(q\in\N, q>1)$. In this construction
Lemma~\ref{l:1/q}
will play a key role as it did in the proof of
Corollary~\ref{cor:euler}.
Therefore we prepare the given polynomials $p_k$ in the following
way: If $P(x)=\sum\limits_{j=0}^N a_j\,x^j$ is any polynomial with
$a_N\neq 0$, then we decompose $P$ by the $q$ polynomials
\[
P_m(x):=\sum_{\begin{array}{c}
\\[-6.5mm]
{\scriptstyle{j=0}}\\[-2.5mm]
{\scriptstyle{j\equiv m \;({\rm mod}\; q)}}
\end{array}}^N
a_j\,x^j
\quad\quad\quad
(m=0,1,\ldots,q-1)
\]
in the form
\[
P(x)=\sum_{m=0}^{q-1} P_m(x)
\;,
\]
and $P$ is the zero polynomial if and only if $P_m$ is the zero
polynomial
for all $m=0,1,\ldots,q-1$. If $m\in\{0,1,\ldots,q-1\}$ then 
$j\equiv m \;({\rm mod}\; q)$ if and only if $j=s_j\,q+m$ with
$s_j\in N_0$.
Therefore we have
\[
P_m(x)=x^m\sum_{\begin{array}{c}
\\[-6.5mm]
{\scriptstyle{j=0}}\\[-2.5mm]
{\scriptstyle{j\equiv m \;({\rm mod}\; q)}}
\end{array}}^N
a_j\,x^{s_j q}
\;.
\]
If we set 
\[
P_m^*(y):=\sum_{\begin{array}{c}
\\[-6.5mm]
{\scriptstyle{j=0}}\\[-2.5mm]
{\scriptstyle{j\equiv m \;({\rm mod}\; q)}}
\end{array}}^N
a_j\,y^{s_j}
\]
we have $P_m(x)=x^m\,P_m^*(x^q)$, and therefore we arrive at the
decomposition $P(x)=\sum\limits_{m=0}^{q-1} x^m\,P_m^*(x^q)$
with the polynomials $P_m^*$. For the degrees of the polynomials
$P_m^*(y)$
we obtain the relations
\be
{\rm deg\:} P_m^*(y)
\left\{
\begin{array}{clc}
\leq &s_N& (0\leq m < r_N)
\\
= &s_N& (m=r_N)
\\
\leq &s_N\!-\!1& (r_N<m\leq q-1)
\end{array}\right.
\;,
\label{eq:degreerelations}
\ee
and $P_{r_N}^*(y)\not\equiv 0$,
where we use the representation
$N=s_N\,q+r_N\;(s_N\in\N_0, r_N\in\N_0, 0\leq r_N\leq q-1)$.

Now we continue with (\ref{eq:givenDE}). We differentiate
(\ref{eq:givenDE})
iteratively $l^*:=n(q-1)$ times to get the equations
\be
\sum_{k=0}^{n+l} p_{kl}(x)\,f^{(k)}(x)=0
\label{eq:first form}
\ee
with polynomials $p_{kl}\;(l=0,1,\ldots,l^*)$,
and we note that $p_{n+l,l}=p_n\not\equiv 0$, and that
$p_{k0}=p_k$.
These form a set of $l^*+1=n(q-1)+1$ equations.

We decompose each polynomial $p_{kl}\;(l=0,1,\ldots, l^*)$ in the
way
introduced above to get
\[
p_{kl}(x)=\sum_{m=0}^{q-1} x^m\,Q_{klm}(x^q)
\]
with polynomials $Q_{klm}$. This brings our $l^*+1$ equations
(\ref{eq:first form}) into the form
\be
\sum_{k=0}^{n+l} \sum_{m=0}^{q-1} Q_{klm}(x^q)\,\lk
x^m\,f^{(k)}(x)\rk=0
\quad\quad\quad(l=0,1,\ldots,l^*)
\;.
\label{eq:second form}
\ee
If we write $k=s_k\,q+r_k\;(s_k\in\N_0,\;r_k=0,1,\ldots,q-1)$,
a brief look in Lemma~\ref{l:1/q} shows that for the purpose to
arrive
at a differential equation for $h$ the ``variables'' to be
eliminated
in the above equations are the terms $\lk x^m\,f^{(k)}(x)\rk$
where $m\in\{0,1,\ldots,q-1\}\setminus \{r_k\}$. Every derivative
$f^{(k)}$
so generates $q-1$ unknowns, and as we have $nq+1$ derivatives,
these are
$(nq+1)(q-1)=q(l^*+1)-1$ unknowns. Therefore we need to generate
$q$ equations out of each of our $l^*+1$ equations (\ref{eq:second
form})
to arrive at a set of $q(l^*+1)$ equations, i.\ e.\ enough
equations to yield the differential equation searched for.

This is done by multiplying each of the equations (\ref{eq:second
form})
by the factors $x^p\;(p=0,\ldots,q-1)$ generating the $q(l^*+1)$
equations
\be
\sum_{k=0}^{n+l} \sum_{m=p}^{q+p-1} Q_{k,l,m-p}(x^q)\,\lk
x^m\,f^{(k)}(x)\rk
=
\sum_{k=0}^{n+l} \sum_{j=0}^{q-1} Q_{klpj}(x^q)\,\lk
x^j\,f^{(k)}(x)\rk
=0
\;,
\label{eq:third form}
\ee
\[
(l=0,1,\ldots,l^*,\;k=s_k\,q+r_k,\;s_k\in\N_0,\;r_k=0,1,\ldots,q-
1,\;
p=0,1,\ldots,q-1)
\]
of the $q(l^*+1)-1$ unknowns $\lk x^j\,f^{(k)}(x)\rk$,
where we use the abbreviations
\be
Q_{klpj}(y):=\funkdeff
{y\,Q_{k,l,q+j-p}(y)}{j=0,\ldots,p-1}
{Q_{k,l,j-p}(y)}{j=p,\ldots,q-1}
\label{eq:abbreviations}
\ee
and
\be
Q_{kl0j}(y):=Q_{klj}(y)
\quad\quad\quad(j=0,1,\ldots,q-1)
\;.
\label{eq:Nebenbed}
\ee
We will show that these unknowns always can be eliminated,
therefore
arriving at a differential equation of order $nq$ for $h$.

We note, that if we differentiate (\ref{eq:givenDE}) iteratively, 
generating $q$ new equations at each step by multiplication with 
$x^p\;(p=0,\ldots,q-1)$, and checking at each step if the unknowns
can be 
eliminated, this algorithm obviously results in the differential
equation of lowest order valid for $h$.

Now we show that in the last step the algorithm always succeeds.
Therefore we rewrite the equations system (\ref{eq:third form})
using the $q\times q$ matrices
\[
A_{kl}(y):=
\Big( Q_{klpj}(y)\Big)_{pj}
\;.
\]
With the aid of these matrices each block of equations with fixed
$l$
can be written as the matrix equation
\[
\sum_{k=0}^{n+l} A_{kl}(x^q)
\lk
\!\!\!
\begin{array}{c}
f^{(k)}(x)\\
x f^{(k)}(x)\\
x^2 f^{(k)}(x)\\
\vdots\\
x^{q-1} f^{(k)}(x)\\
\end{array}
\!\!\!
\rk=0
\;.
\]
Because of (\ref{eq:abbreviations}) and (\ref{eq:Nebenbed}) the
matrices
$A_{kl}(y)$ are of the following special form
\[
A_{kl}(y)=
\lk
\begin{array}{cccccc}
Q_{kl0}(y)&Q_{kl1}(y)&Q_{kl2}(y)&\cdots&Q_{k,l,q-2}(y)&Q_{k,l,q-1
}(y)\\
y\,Q_{k,l,q-1}(y)&Q_{kl0}(y)&Q_{kl1}(y)&\cdots&Q_{k,l,q-3}(y)&Q_{
k,l,q-2}(y)\\
y\,Q_{k,l,q-2}(y)&y\,Q_{k,l,q-1}(y)&Q_{kl0}(y)&\cdots&Q_{k,l,q-4}
(y)&
Q_{k,l,q-3}(y)\\
\vdots&\vdots&\vdots&\ddots&\vdots&\vdots\\
y\,Q_{kl2}(y)&y\,Q_{kl3}(y)&y\,Q_{kl4}(y)&\cdots&Q_{kl0}(y)&Q_{kl
1}(y)\\
y\,Q_{kl1}(y)&y\,Q_{kl2}(y)&y\,Q_{kl3}(y)&\cdots&y\,Q_{k,l,q-1}(y
)&Q_{kl0}(y)
\end{array}
\rk
\;.
\]
From the relations $\sum\limits_{k=0}^{n+l}
p_{kl}(x)\,f^{(k)}(x)=0$ with
$p_{kl}(x)=x^m \sum\limits_{m=0}^{q-1}Q_{klm}(x^q)$, 
$p_{n+l,l}(x)=p_n(x)\not\equiv 0$, and $p_{k0}(x)=p_k(x)$ we see
that
the matrix $A_{kl}(y)$ is completely determined by the polynomial
$p_{kl}(x)$,
and we have
\[
A_{n+l,l}(y)=A_{n0}(y)=:A_n(y)
\;.
\]
In particular, the matrix $A_n(y)$ is completely determined by the
nonvanishing
coefficient polynomial $p_n(x)$ of $f^{(n)}$ in the original
equation
(\ref{eq:givenDE}). We will show now that the matrix $A_n(y)$ is
invertible.
Let $p_n(x)=\sum\limits_{j=0}^N a_j\,x^j$ with $a_N\neq
0,\;N=s_N\,q+r_N\;
(s_N, r_N\in\N_0,\;0\leq r_N\leq q-1)$. We have 
$p_n(x)=p_{n0}(x)=x^m\sum\limits_{m=0}^{q-1}Q_{n0m}(x^q)$. If we
use
the notation $Q_m(y):=Q_{n0m}(y)$, then for $Q_m$ the degree
relations 
(\ref{eq:degreerelations}) read
\[
{\rm deg\:} Q_m(y)
\left\{
\begin{array}{clc}
\leq &s_N& (0\leq m < r_N)
\\
= &s_N& (m=r_N)
\\
\leq &s_N\!-\!1& (r_N<m\leq q-1)
\end{array}\right.
\;,
\]
and we have $Q_{r_N}\not\equiv 0$. This information leads to
the following observations
about the degrees of the polynomial entries of the matrix $A_n$:
The
entries of $A_n$ for which the difference of colomn index and
row index $\Delta$ is $\geq (r_N+1)$ (region I) have degrees $\leq
s_N-1$, 
the entries for which $0\leq\Delta\leq r_N$ (region II) have
degrees $\leq s_N$,
further if $-r_N+1\leq \Delta\leq -1$ (region III), then the
degrees again are 
$\leq s_N$, and finally if $\Delta\leq -r_N$ (region IV), then
the degrees are $\leq s_N+1$.

Assume next that $p_n(0)\neq 0$, then we get $|A_n(0)|=Q_0^q(0)=
p_n^q(0)\neq 0$, and therefore $A_n(y)$ is invertible. To obtain
the same result if $p_n(0)=0$, we use our degree observations. By
the
Weierstra{\ss} representation of the determinant $D=|d_{jk}|$ of a
$q\times q$  matrix $(d_{jk})$
\[
D=\sum_{\pi\in S_q} {\rm sign\:} \pi \,
d_{1,\pi(1)}\,d_{2,\pi(2)}\cdots \,d_{q,\pi(q)}
\]
where the sum is to be taken over all permutations $\pi\in S_q$,
we observe that in the Weierstra{\ss} representation of the 
determinant of our matrix $A_n(y)$ the summands
\[
Q_0^q(y),\pm y\,Q_1^q(y)
,\ldots,
\pm y^{r_N}\,Q_{r_N}^q(y),\ldots,
y^{q-1}\,Q_{q-1}^q(y)
\]
occur. Now it follows that
\[
{\rm deg\:}\lk y^{r_N}\,Q_{r_N}^q(y)\rk=q\,s_N+r_N=N={\rm deg\:}
p_n(x)>0
\;,
\]
and for all other summands $P(y)$ 
using our degree observations we get the relation
\[
{\rm deg\:} P(y) <q\,s_N+r_N=N
\;.
\]
Hence the degree of $|A_n(y)|$ equals $N$, thus $|A_n(y)|$ cannot
be the 
zero polynomial, and $A_n(y)$ therefore is invertible.

Using the inverse $A_n^{-1}$, we can now write
\[
\lk
\begin{array}{c}
f^{(n+l)}(x)\\
x f^{(n+l)}(x)\\
x^2 f^{(n+l)}(x)\\
\vdots\\
x^{q-1} f^{(n+l)}(x)\\
\end{array}
\!\!\!
\rk=
-\sum_{k=0}^{n+l-1} A_n^{-1}(x^q)\cdot A_{kl}(x^q)
\lk
\begin{array}{c}
f^{(k)}(x)\\
x f^{(k)}(x)\\
x^2 f^{(k)}(x)\\
\vdots\\
x^{q-1} f^{(k)}(x)\\
\end{array}
\!\!\!
\rk
\quad\quad(l=0,\ldots,l^*)
\;.
\]
Using the fact that the inhomogeneous parts of these equations do
not vanish,
it is now easily established that the unknowns $\lk
x^j\,f^{(k)}(x)\rk$
can be eliminated, which finishes the algorithm.\hfill$\Box$
\\[3mm]
After these preparations we are ready to state our main theorem
(compare \cite{Sta}, \cite{SZ}).
\bT
\label{th:main theorem}
{\rm
Let the two functions $f$ and $g$ be simple of order $n$ and $m$,
respectively, and let $r$ be rational. 
Then 
\one
{a}
{\int f(x)\,dx\;\;\;\mbox{is simple of order $\leq n+1$}\;,}
\one
{b}
{f'\;\;\;\mbox{is simple of order $\leq n$}\;,}
\one
{c}
{f+g\;\;\;\mbox{is simple of order $\leq n+m$}\;,}
\one
{d}
{f\cdot g\;\;\;\mbox{is simple of order $\leq n\,m$}\;,}
\one
{e}
{f\circ r\;\;\;\mbox{is simple of order $\leq n$}\;,}
\one
{f}
{f\circ x^{p/q}\quad(p,q\in\Z)\;\;\;\mbox{is simple of order $\leq
n\,q$}\;.}
}
\eT
\pr
(a): Let $h=\int f(x)\,dx$, and let $f$ satisfy the simple
differential
equation (\ref{eq:simple DE}), then obviously
\[
\sum\limits_{k=0}^n p_k\,h^{(k+1)}=0
\;,
\]
and so $h$ is simple of order $\leq n+1$.
\\
(b):
Let $h=f'$, and (\ref{eq:simple DE}) be valid. If $p_0\equiv 0$
then
(\ref{eq:simple DE}) is a simple differential equation valid for
$h$, and so $h$ is simple of order $n$. 
If $p_0\not\equiv 0$, however, then, dividing by $p_0$,
we write (\ref{eq:simple DE}) as
\be
f=\sum\limits_{k=1}^n r_k\,f^{(k)}
\label{eq:fsub}
\ee
with rational functions $r_k$. If we differentiate the original
differential
equation, we get a further differential equation
\[
\sum\limits_{k=0}^{n+1} P_k\,f^{(k)}=0
\]
with polynomials $P_k\;(k=0,\ldots,n+1)$. We substitute $f$ by
(\ref{eq:fsub}),
and multiply by the common denominator to get a simple differential
equation of order $\leq n$ for $h$. 
\\
(c):
First we give an algebraic argument for the existence statement.
Let $f$ and $g$ satisfy simple differential equations of order $n$
and
$m$, respectively.  
We consider the linear space $L_f$ of functions with rational
coefficients
generated by $f, f',f'',\ldots,f^{(k)},\ldots$. As $f,
f',\ldots,f^{(n)}$ are linearly dependent by (\ref{eq:simple DE}), 
and as the same conclusion follows for $f',f'',\ldots,f^{(n+1)}$ by
differentiation, and so on inductively,
the dimension of $L_f$ is $\leq n$. Similarly $L_g$ has dimension
$\leq m$. 
We now build the sum $L_f +L_g$ which is of dimension $\leq n+m$.
As $f\!+\!g,(f\!+g)',\ldots,(f\!+\!g)^{(k)},\ldots$ are
elements of $L_f +L_g$, arbitrary $n+m+1$ many of them are linearly
dependent.
In particular, $f+g$ is simple of order $\leq n+m$. 

Now we like to present an algorithm which generates the
differential
equation for the sum: 
We may bring the given differential equations for $f$ and $g$ into
the form
\be
f^{(n)}=\sum_{j=0}^{n-1} p_j\,f^{(j)}
\label{eq:fde}
\;,
\ee
and
\be
g^{(m)}=\sum_{k=0}^{m-1} q_k\, g^{(k)}
\label{eq:gde}
\;,
\ee
with rational functions $p_j$, and $q_k$ respectively. 
By iterative differentiation and recursive substitution
of (\ref{eq:fde}), and (\ref{eq:gde}), we generate sets of rules
\be
f^{(l)}=\sum_{j=0}^{n-1} p_j^l\,f^{(j)}
\quad\quad\quad(l=n,\ldots,n+m)
\label{eq:frules}
\;,
\ee
and
\be
g^{(l)}=\sum_{k=0}^{m-1} q_k^l\, g^{(k)}
\quad\quad\quad(l=m,\ldots,n+m)
\label{eq:grules}
\;,
\ee
with rational functions $p_j^l$, and $q_k^l$, respectively. 
Now, by iterative differentiation of the defining equation
$h:=f+g$, 
we generate the set of equations
\be
\begin{array}{rcl}
h&=&f+g\\
h'&=&f'+g'\\
h''&=&f''+g''\\
\vdots\;\;\;&&\;\;\;\vdots\\
h^{(n+m)}&=&f^{(n+m)}+g^{(n+m)}
\end{array}
\;.
\label{eq:hequations}
\ee
Next, we take the first $\max\{n,m\}$ of these equations, and use
the rules
(\ref{eq:frules}) and (\ref{eq:grules}) to eliminate all
occurrences of
$f^{(l)}\;(l\geq n)$, and $g^{(l)}\;(l\geq m)$. On the right hand
sides
remain the $n+m$ variables 
$f^{(l)}\;(l=0,\ldots,n-1)$ and $g^{(l)}\;(l=0,\ldots,m-1)$.
We solve the remaining linear system of equations for the variables
$h^{(l)}\;(l=0,\ldots,\max\{n,m\})$ trying to eliminate the
variables
$f^{(l)}\;(l=0,\ldots,n-1)$ and $g^{(l)}\;(l=0,\ldots,m-1)$.
If this procedure is successful, it generates the simple
differential
equation for $h$, searched for.

If the procedure fails, however,
we go on taking the first $\max\{n,m\}+1$ of equations 
(\ref{eq:hequations}), doing the same manipulations,
and so on, until we take the whole set of equations
(\ref{eq:hequations}), where the procedure must stop.
\\
(d):
This situation is handled in exactly the same way as case (c). The
only 
difference is that to differentiate the product $h:=f g$ we use the
Leibniz rule, and the variables to be eliminated are the $n\,m$
products $f^{(j)}\,g^{(k)}\;(j=0,\ldots,n-1,\;k=0,\ldots,m-1)$. It
is
easily seen that this procedure
generates a differential equation for $h$ of order $\leq n\,m$.

That the algorithm stops is seen by the fact that finally we arrive
at
a set of $n\,m+1$ equations with $n\,m$ variables to be eliminated,
where all equations actually possess a nonvanishing
inhomogeneous part $h^{(l)}$.
\\
(e):
Let $h:=f\circ r$ for some rational $r$, and let $f$ satisfy
(\ref{eq:simple DE}). We compose (\ref{eq:simple DE}) with $r$, and
get
\be
0=\sum\limits_{k=0}^n (p_k\!\circ\! r)\cdot (f^{(k)}\!\circ\! r)
=\sum\limits_{k=0}^n q_k\cdot \lk f^{(k)}\!\circ\! r\rk
\quad\quad(q_k\;\;\mbox{rational functions})
\;.
\label{eq:DE subst r}
\ee
Differentiating the identity $h=f\circ r$ leads to
\beao
h'=f'\!\circ\! r\cdot r'
&\quad\mbox{or}\quad&
f'\!\circ\! r=h'/r'
\;,
\\
h''=f'\!\circ\! r\cdot r''+f''\!\circ\! r\cdot r'^2
&\quad\mbox{or}\quad&
f''\!\circ\! r=\ed{r'^2}\lk h''-f'\!\circ\! r\cdot r''\rk
=
\ed{r'^2}\lk h''-\frac{r''}{r'}h'\rk
\;,
\eeao
and inductively it follows
\[
h^{(k)}=\sum_{j=1}^k r_j\cdot \lk f^{(j)}\!\circ\! r\rk
\quad\quad(r_j\;\;\mbox{rational functions})
\;.
\]
Solving for $f^{(k)}\!\circ\! r$ and substituting the previous
results we 
are lead to a representation
\[
f^{(k)}\!\circ\! r=
\sum_{j=1}^k R_j\cdot h^{(j)}
\]
with rational functions $R_j$.

Substituting these results in (\ref{eq:DE subst r}) and multiplying
by
the common denominator yields a simple differential equation for
$h$,
and we see that the order of $h$ is $\leq n$.
\\
(f):
The statement follows immediately by an application of
Theorem~\ref{th:1/q}
and (e).
\eop
\begin{remark}
{\rm
We remark that again our proofs provide algorithms to generate the
differential equations searched for.

Further we note that by Theorem~\ref{th:tanx} the functions
$1/g$ and $f/g$ in general are not simple, if $f$ and $g$
are, as $\sec x=1/\cos x$ and $\tan x=\sin x/\cos x$ show.

Note that Theorem~\ref{th:main theorem} (f) generalizes Theorem 8.1
in
\cite{Koe92}.
}
\end{remark}
The following consideration strengthens Theorem~\ref{th:main
theorem} (c),
and brings it in connection with the fundamental systems of simple
differential equations.
\bC
{\rm
Assume, $f$ satisfies a simple differential equation $F$ of order
$n$,
$g$ satisfies a simple differential equation $G$ of order $m$, and
$F$ and $G$
do not have any {\sl common nontrivial} solution. Then the order of
any
simple differential equation satisfied by the sum $h:=f+g$ is $\geq
n+m$.
In particular, $h$ satisfies no simple differential equation of
order 
$<n+m$, and exactly one differential equation $H$ of order $n+m$
which is 
generated by the algorithm given in Theorem~\ref{th:main theorem}
(c).
}
\eC
\pr
Let $f_1, f_2,\ldots, f_n$ be a fundamental system of the
differential 
equation $F$, and $g_1, g_2, \ldots, g_m$ be a fundamental system
of
the differential equation $G$. Then $f_1, f_2,\ldots, f_n$ are
linearly independent (in the usual sense), and $g_1, g_2, \ldots,
g_m$ are
linearly independent, as well. Assume now, the $n+m$ functions
$f_1, f_2,\ldots, f_n, g_1, g_2, \ldots, g_m$ are linearly
dependent.
Then there are constants $a_j\;(k=j,\ldots,n)$ and
$b_k\;(k=1,\ldots,m)$,
not all equal to zero, such that
\be
\sum_{j=1}^n a_j\,f_j+\sum_{k=1}^m b_k\,g_k=0
\;.
\label{eq:agree}
\ee
From the linear independence of the subsets $f_1, f_2,\ldots, f_n$,
and
$g_1, g_2, \ldots, g_m$ it follows that at least one of the numbers

$a_j\;(j=1,\ldots,n)$, and one of the numbers $b_k\;(k=1,\ldots,m)$
does
not vanish. Therefore by the linearity of $F$ and $G$ the functions
$f:=\sum\limits_{j=1}^n a_j\,f_j$ and $g:=-\sum\limits_{k=1}^m
b_k\,g_k$
are solutions of $F$, and $G$, respectively,
which by (\ref{eq:agree}) agree, and we have a contradiction.

Thus we have proved that $f_1, f_2,\ldots, f_n, g_1, g_2, \ldots,
g_m$ are
linearly independent. As the zero function is a solution of both
$F$, and $G$,
it turns out that $f_1, f_2,\ldots, f_n, g_1, g_2, \ldots, g_m$
must satisfy
any linear differential equation $H$ for the sum, and so form a
fundamental
system for $H$. Therefore the order of $H$ is $\geq n+m$.
\eop
We mention that, on the other hand, if the algorithm presented in
Theorem~\ref{th:main theorem} (c) generates a differential equation
of order $<m+n$ for $h:=f+g$, then the 
simple differential equations $F$ for $f$ and $G$ for $g$ must have
a common nonzero solution, i.\ e.\ the fundamental systems $F$ and
$G$
are linearly dependent.
We gain this insight without solving any of the differential
equations!

From Theorem~\ref{th:main theorem} (d) a similar statement can be
obtained 
for the product $h:=f g$: If our algorithm generates a differential
equation of order $<n\,m$ for $h$, then the product of the
fundamental systems $S_F$ and $S_G$ of $F$ and $G$, i.\ e.\ the set

$\{f\,g\;|f\in S_F,\;g\in S_G\}$, is linearly dependent.
Again, we gain this insight without solving any of the differential
equations.

In connection with Theorem~\ref{th:main theorem} (d) for the
product $h:=f g$ 
we note that for the special case of the square $h:=f^2=f\cdot f$ 
the bound for the order of $h$ can be considerably strengthened.
\bC
\label{cor:square}
{\rm
Let $f$ be simple of order $n$, then $h:=f^2$ is simple of order
$\leq
\frac{n(n+1)}{2}$.
}
\eC
\pr
A careful study of the algorithm given in the proof of
Theorem~\ref{th:main theorem} (d) shows that the order cannot
exceed 
$\frac{n(n+1)}{2}$. This depends on the fact that
the variables to be eliminated are the $\frac{n(n+1)}{2}$ different
products $f^{(j)}\,f^{(k)}\;(j=0,\ldots,n-1,\;k=0,\ldots,n-1)$. 
\eop

\section{Algorithmic calculation of simple differential equations}

Theorem~\ref{th:main theorem} enables us to define the linear space
$F$ of
functions 
generated by the algebraic functions, $\exp x$, $\ln x$, $\sin x$,
$\cos x$, $\arcsin x$, $\arctan x$, (and, if we wish, further
special functions satisfying simple differential equations,)
and the functions that are constructed by an application of
finitely many of the procedures
(a)--(f) of Theorem~\ref{th:main theorem}. The theorem then states
in
particular that $F$ forms a differential ring.

The algorithm which generates the simple differential equation of
lowest
order for $f$ given in \cite{Koe92}--\cite{Koe93c} is strengthened
by
Theorem~\ref{th:main theorem}. 
\begin{algorithm}[Find a simple differential equation]
\label{algorithm:Find a simple DE}%
{\rm
Let $f\in F$. Then the following procedure generates the simple
differential equation of lowest order valid for $f$:
\begin{enumerate}
\item[{\rm (a)}]
Find out whether there exists a simple differential equation 
for $f$ of order $N:=1$.
Therefore differentiate $f$, and solve the linear equation
\[
f'(x)+A_{0}f(x)=0
\]
for $A_{0}$; i.\ e.\ set $A_{0}:=-\frac{f'(x)}{f(x)}$.
Is $A_{0}$ rational in $x$, then you are done after
multiplication with its denominator.
\item[{\rm (b)}]
Increase the order $N$ of the differential equation searched for by
one.
Expand the expression
\[
f^{(N)}(x)+A_{N-1}f^{(N-1)}(x)+\cdots+A_{0}f(x)
\;,
\]
and check, if the summands contain exactly $N$
rationally independent expressions (i.\ e.: linearly independent
over $K[x]$)
considering the numbers $A_{0}, A_{1},\ldots, A_{N-1}$ as
constants.
Just in that case there exists a solution as follows: Sort with
respect
to the rationally independent terms and create a system of linear
equations
by setting their coefficients to zero. Solve this system for the
numbers
$A_0, A_1,\ldots, A_{N-1}$. Those are rational functions in $x$,
and
there
exists a unique solution. After multiplication by the common
denominator
of
$A_0, A_1,\ldots, A_{N-1}$ you get the differential equation
searched for.
Finally cancel common factors of the polynomial coefficients.
\item[{\rm (c)}]
If part (b) was not successful, repeat step (b).
\end{enumerate}
}
\end{algorithm}
\pr
The proof given in \cite{Koe92} shows that the given algorithm
indeed
generates the simple differential equation of lowest order valid
for $f$
whenever such a differential equation exists. Theorem~\ref{th:main
theorem} now guarantees the existence of such a differential
equation, 
and therefore the algorithm does not end in an infinite loop.
\eop
\begin{remark}
{\rm
Obviously $f\in F$ can be checked by a pattern matching mechanism
applied to the given expression $f$, and
Theorem~\ref{th:main theorem} then moreover gives a priori
estimates for the 
possible order of the solution differential equation.

We further point out that an actual implementation of the algorithm
should be able to decide the rational dependency of a set of
functions.
Otherwise, it may fail to find the simple differential equation of
lowest 
order. An example of this type is given by
$f(x):=\sin\:(2x)-2\,\sin x\cos\ x$,
for which the algorithm yields the differential equation $f''+f=0$
rather
than $f'=0$ if the rational dependency of $\sin\:(2x)$, and $\sin
x\cos x$
is not discovered.
}
\end{remark}
In a forthcoming paper \cite{Koe94} we discuss the more general
situation of
functions $f$ that depend on other special functions not contained
in $F$,
like Airy functions, Bessel functions and orthogonal polynomials.

We mention that the algorithms presented in Theorem~\ref{th:main
theorem}
can be combined to get the following different

\begin{algorithm}[Find a simple differential equation]
\label{algorithm:Find a simple DE2}%
{\rm
Let be $f\in F$. Then recursive application of the algorithms
presented in 
Theorem~\ref{th:main theorem} generates a simple differential
equation 
for $f$. 
}
\end{algorithm}
\begin{remark}
{\rm
We note that, in general, we cannot control the order of the
resulting
differential equation as the algorithm uses recursive descent
through 
the expression tree of the given $f$. Given a sum of two functions
or order $n$ and $m$, however,
it is easily seen that the order of the resulting differential
equation
is the lowest possible order which is $\geq \max\{n,m\}$. So with
the given
algorithm, it is principally impossible to derive any differential
equation 
of order lower than two for the example function 
$f(x):=\sin\:(2x)-2\,\sin x\cos\ x$ as the lowest order
differential equations
for the summand functions $\sin\:(2x)$, and $\sin x\cos\ x$ are of
order two.

This can be interpreted as follows: By construction, for a sum
$f:=f_1+f_2$
a differential equation is obtained that is valid not only for
$f_1+f_2$,
but for any linear combination $\la_1\,f_1+\la_2\,f_2$ or $f_1$ and
$f_2$.
In other words, this subalgorithm generates a differential equation
for the 
linear hull of $f_1$ and~$f_2$.
}
\end{remark}

through the
quotient

\section{Sharpness of the orders}

In this section we give examples that show that the bounds for the
least 
orders that are given in Theorems~\ref{th:1/q}--\ref{th:main
theorem},
and in Corollary~\ref{cor:square} in fact may be assumed.

First we consider the statement of Corollary~\ref{cor:square}. We
get
\bT
\label{th:sharp square}
{\rm
For each $n\in\N$ there is a function $f_n$ that is simple of 
order $n$, such that $h:=f_n^2$ is simple of order
$\frac{n(n+1)}{2}$.
Indeed, an example function of that type is
\be
f_n(x):=\sum\limits_{k=1}^{n} e^{x^k}
\;.
\label{eq:squareexample}
\ee
}
\eT
\pr
We show that $f_n$, given by (\ref{eq:squareexample}), is simple of
order $n$, and that $f_n^2$ is simple of order $\frac{n(n+1)}{2}$.

As each function $e^{x^k}\;(k=1,\ldots,n)$ is simple of order $1$
(satisfying the differential equation $f'-k\,x^{k-1}\,f=0$), 
by an inductive application of Theorem~\ref{th:main theorem} (c),
we
see that $f_n$ is simple of order $\leq n$. To show that $f_n$ is
simple of 
order $\geq n$, we first show that $e^{x^k}\;(k=1,\ldots,n)$ are
rationally
independent. Suppose, any linear combination 
\be
r_1\,e^{x}+r_2\,e^{x^2}+\cdots + r_{n}\,e^{x^{n}}=0
\label{eq:linear combination}
\ee
with rational functions $r_k\;(k=1,\ldots,n)$
representing zero is given. From (\ref{eq:linear combination}) we
get,
assuming $r_k\not\equiv 0$
\[
r_k\,e^{x^k}=-\sum_{\begin{array}{c}
\\[-6.5mm]
{\scriptstyle{j=1}}\\[-2.5mm]
{\scriptstyle{j\neq k}}
\end{array}}^{n}
r_j\,e^{x^j}
\;,
\]
and dividing by $r_k\,e^{x^k}$, we get
\[
1=-\sum_{\begin{array}{c}
\\[-6.5mm]
{\scriptstyle{j=1}}\\[-2.5mm]
{\scriptstyle{j\neq k}}
\end{array}}^{n}
R_j\,e^{x^j-x^k}
\]
with rational functions $R_j\;(j=1,\ldots,n)$.
For $k=n$, we get by taking the limit $x\pf\infty$, along the
positive
real axis
\[
1=-\lim_{x\pf\infty} \sum_{j=1}^{n-1} R_j\,e^{x^j-x^{n}}=0
\;,
\]
a contradiction to the assumption $r_{n}\not\equiv 0$.
For $k<n$, the right hand limit
\[
-\lim_{x\pf\infty}\sum_{\begin{array}{c}
\\[-6.5mm]
{\scriptstyle{j=1}}\\[-2.5mm]
{\scriptstyle{j\neq k}}
\end{array}}^{n}
\,R_j\,e^{x^j-x^k}
\]
does not exist, whereas the left hand side equals 1,
a contradiction to the assumption $r_k\not\equiv 0$. Therefore 
$r_k\equiv0\;(k=1,\ldots,n)$, and $e^{x^k}\;(k=1,\ldots,n)$ are
rationally
independent.

From this it follows that applying Algorithm~\ref{algorithm:Find a
simple DE}
to $f_n$, in each step $n$ rationally independent are involved (the
derivative
of $e^{x^k}$ is of the same type), so that the algorithm can be
successful
not earlier than in the $n$th step. As this algorithm always finds
the simple
differential equaton of least order, the order of $f_n$ is $\geq
n$.

Next we consider the function $f_n^2$. Expanding $f_n^2$ we get
the representation
\[
f_n^2(x)=
e^{2x}+e^{2x^2}+\cdots+e^{2x^{n}}+2\,e^x\,e^{x^2}+\cdots+
2\,e^{x^{n-1}}\,e^{x^{n}}
\]
consisting of $n+\ueber{n}{2}=\frac{n(n+1)}{2}$ summands. An
application
of Theorem~\ref{th:main theorem} (c), again, shows that $f_n^2$ is
simple of order $\leq \frac{n(n+1)}{2}$, as each summand 
$f:=e^{x^j}\,e^{x^k}$ satisfies the simple differential equation
\[
f'-2\,(j\,x^j+k\,x^k)\,f=0
\]
of order $1$.
Similarly as above one realizes
that $e^{x^j}\,e^{x^k}\;(j,k=1,\ldots,n \;(j\neq k))$
are rationally independent, and therefore by an application of
Algorithm~\ref{algorithm:Find a simple DE} we conclude that $f_n^2$
is
simple of order $\geq \frac{n(n+1)}{2}$, finishing the proof.
\eop
Next we consider the statement of Theorem~\ref{th:1/q}. We get
\bT
\label{th:sharp 1/q}
{\rm
For each $n,q\in\N$ there is a function $f_{nq}$ that is simple of
order $n$, such that $h:=f\circ x^{1/q}$ is simple of order $n\,q$.
Indeed, an example function of that type is
\be
f_{nq}(x):=\sum\limits_{k=1}^{n} e^{x^{kq+1}}
\;.
\label{eq:1/qsharp}
\ee
}
\eT
\pr
Everything follows as in Theorem~\ref{th:sharp square}, if we show
that 
\[
f(x):=e^{x^{kq+1}}\quad\quad (k=1,\ldots,n,\; q\in\N)
\]
is simple of order $1$, whereas
\[
h(x):=e^{x^{k+1/q}}\quad\quad (k=1,\ldots,n,\; q\in\N)
\]
is simple of order $q$.
The first statement was already proven so that it remains to prove
the
second statement. This, however, turns out to be
a consequence, of Lemma~\ref{l:1/q} (a).
We start with one rationally independent function, namely $h$
itself.
As $h(x)=f(x^{1/q})$, Lemma~\ref{l:1/q} (a) shows that in the $r$th
differentiation step ($r\leq q$) for $h$ the new summand
\[
x^{r/q}\,f^{(r)}(x^{1/q})
\]
is introduced. 
By a similar argument as in Theorem~\ref{th:sharp square} one sees
that
for each $r\leq q-1$ the functions involved
functions
\[
\left\{
h(x),
x^{1/q}\,f'(x^{1/q}),
\ldots,
x^{r/q}\,f^{(r)}(x^{1/q})
\right\}
\]
are rationally independent.

In the $q$th step, however, by (\ref{eq:scharf1/q}),
$h^{(q)}$ rationally depends on $h, h',\ldots, h^{(q-1)}$,
so that an arbitrary
linear combination of $h, h',\ldots h^{(q)}$ contains exactly $q$
rationally
independent terms, showing that $h$ is of order $q$.
\eop
Finally we consider the statements of Theorem~\ref{th:main
theorem}.
We get
\bT
\label{th:sharp main (c)}
{\rm
The statements of Theorem~\ref{th:main theorem} 
are sharp which is seen by the following example functions:
\begin{enumerate}
\item[(a)]
For each $n$ the function $f_n$ given by (\ref{eq:squareexample})
is simple
of order $n$, and $\int f_n$ is simple of order $n+1$.
\item[(b)]
For each $n$ the function $f_n$ given by (\ref{eq:squareexample})
is simple
of order $n$, and $f_n'$ is simple of order $n$, too.
\item[(c)]
For each $n,m\in\N$ there are functions $f$ and $g$ that are simple
of
order $n$, and $m$, respectively,
such that $h:=f+g$ is simple of order $n+m$.
Indeed, example functions of that type are
\be
f(x):=\sum\limits_{k=1}^{n} e^{x^{k}}
\quad\quad\mbox{and}\quad\quad
g(x):=\sum\limits_{k=n+1}^{n+m} e^{x^{k}}
\;.
\label{eq:f+gexample}
\ee
\item[(d)]
For each $n,m\in\N$ there are functions $f$ and $g$ that are simple
of
order $n$, and $m$, respectively,
such that $h:=f\, g$ is simple of order $n\, m$.
Indeed, example functions of that type are given by
(\ref{eq:f+gexample}).
\item[(e)]
For each $n$ the function $f_n$ given by (\ref{eq:squareexample})
is simple
of order $n$, and for each rational function $r$ the function
$f\circ r$ is simple of order $n$, too.
\item[(f)]
For each $n$ the function $f_{nq}$ given by (\ref{eq:1/qsharp}) is
simple
of order $n$, and for each $p\in\N$ the function
$f\circ x^{p/q}$ is simple of order $n\,q$.
\end{enumerate}
}
\eT
\pr
(a) By Theorem \ref{th:sharp square} $f_n$ is simple of order $n$.
As the antiderivative of $f_n$ is not an elementary function, it is
rationally independent, and so its order is $\geq n+1$.
\\
(b)
As the derivatives of each summand of $f_n$ depend rationally on
$f_n$,
by the linearity the order of $f_n'$ equals $n$, too.
\\
(c)
In Theorem~\ref{th:sharp square} we already proved that $f$ is
simple
of order $n$, and that $f+g$ is simple of order $n+m$. 
A similar argument shows that $g$ is simple of order $m$, and we
are done.
\\
(d)
If we expand the function $f\,g$, we get
\[
f(x)\,g(x)=\sum\limits_{k=1}^{n} e^{x^{k}}\,
\sum\limits_{k=n+1}^{n+m} e^{x^{k}}=
\sum\limits_{k=1}^{n}\sum\limits_{j=1}^m e^{x^{k}}\,e^{x^{j+n}}
\]
with $n\,m$ rationally independent expressions. Thus the order of
$f\,g$ is
$\geq n\,m$.
\\
(e)
By the argumentation of Theorem~\ref{th:sharp square} one sees that
$f_n\circ r$ has order $n$.
\\
(f)
This follows from (e) and Theorem~\ref{th:sharp 1/q}.
\eop

\end{document}